\documentclass[11pt]{article}
\usepackage{latexsym}

\newcommand{\eproof}{\mbox{\ }\hfill $\Box$ \par \vskip 10pt}

\newtheorem{Theorem}{Theorem}[section]
\newtheorem{lemma}[Theorem]{Lemma}
\newtheorem{prop}[Theorem]{Proposition}

\newtheorem{corol}[Theorem]{Corollary}

\begin{document}

\title{High frequency resolvent estimates for perturbations by large  
long-range magnetic potentials and applications to dispersive estimates}

\author{{\sc Fernando Cardoso, Claudio Cuevas and Georgi Vodev\thanks{Corresponding author}}}

\date{}

\maketitle

\noindent
{\bf Abstract.} We prove optimal high-frequency resolvent estimates for self-adjoint operators of the form
$G=\left(i\nabla+b(x)\right)^2+V(x)$ on $L^2({\bf R}^n)$, $n\ge 3$, where the magnetic potential $b(x)$ and the electric
potential $V(x)$ are long-range and large. As an application, we prove dispersive estimates for the wave group
$e^{it\sqrt{G}}$ in the case $n=3$ for potentials $b(x), V(x)=O(|x|^{-2-\delta})$ for $|x|\gg 1$, where $\delta>0$. 

\setcounter{section}{0}
\section{Introduction and statement of results}

The purpose of the present paper is to study the high frequency behavior of the resolvent of self-adjoint operators on 
$L^2({\bf R}^n)$, $n\ge 3$, of the form
$$G=\left(i\nabla+b(x)\right)^2+V(x),$$
where $b(x)=\left(b_1(x),...,b_n(x)\right)$ is a vector-valued magnetic potential and $V$ is an electric potential, 
$b_j$ and $V$ being real-valued functions. To describe the class these functions 
belong to, we introduce the polar coordinates $r=|x|$, $w=\frac{x}{|x|}\in {\bf S}^{n-1}$. 
They are of the form $b(x)=b^L(x)+b^S(x)$,
$V(x)=V^L(x)+V^S(x)$, where $b^L$ and $V^L$ are $C^1({\bf R}^+)$, ${\bf R}^+=(0,+\infty)$, 
functions with respect to the radial variable $r$. 
We suppose that
there exist constants $C>0$, $0<\delta\ll 1$ so that for all $(r,w)\in {\bf R}^+\times{\bf S}^{n-1}$ we have
$$\left|V^L(rw)\right|\le C,\eqno{(1.1)}$$
 $$\partial_rV^L(rw)\le C\psi_\delta(r),\eqno{(1.2)}$$
$$\left|V^S(rw)\right|\le C\langle r\rangle^{-1-\delta},\eqno{(1.3)}$$
$$\left|\partial_r^kb^L(rw)\right|\le Cr^{1-k}\psi_\delta(r),\quad k=0,1,\eqno{(1.4)}$$
$$\left|b^S(rw)\right|\le C\eta_\delta(r),\eqno{(1.5)}$$
where $\psi_\delta(r)=r^{-1+\delta}\langle r\rangle^{-2\delta}$, $\eta_\delta(r)=r^{\delta}\langle r\rangle^{-1-2\delta}$. 
Finally, we suppose that the function $b^S(rw)$ is continuous in $r$ uniformly in $w$. More precisely, we assume that the function 
$g_\delta(r,w)=b^S(rw)/\eta_\delta(r)$ satisfies
$$\forall\epsilon>0\,\exists\theta=\theta(\epsilon)>0 \,\,\mbox{so that} \,\,
\left|g_\delta(r+\theta\sigma,w)-g_\delta(r,w)\right|
\le \epsilon$$ $$\mbox{for all}\,\, r>0, 0<\sigma\le 1,\,w\in{\bf S}^{n-1}.\eqno{(1.6)}$$
Our main result is the following

\begin{Theorem}Under the assumptions (1.1)-(1.6), for every $\delta'>0$ there exist constants $C,\lambda_0>0$ so that for 
$\lambda\ge\lambda_0$, $0<\varepsilon\le 1$, $0\le|\alpha_1|,|\alpha_2|\le 1$, we have the estimate
$$\left\|\langle x\rangle^{-\frac{1+\delta'}{2}}\partial_x^{\alpha_1}\left(G-\lambda^2\pm i\varepsilon\right)^{-1}\partial_x^{\alpha_2}
\langle x\rangle^{-\frac{1+\delta'}{2}}\right\|_{L^2\to L^2}\le C\lambda^{|\alpha_1|+|\alpha_2|-1}.\eqno{(1.7)}$$
Moreover, if in addition we suppose that $b^S\equiv 0$ and the functions $b=b^L$ and $V=V^L+V^S$ satisfy
$$\left|\frac{\partial(r^2V^L(rw))}{\partial r}\right|\le Cr\psi_\delta(r),\eqno{(1.8)}$$
$$\left|V^S(rw)\right|\le C\langle r\rangle^{-2-\delta},\eqno{(1.9)}$$
$$\left|\partial_r^kb(rw)\right|\le Cr^{-k}\psi_\delta(r),\quad k=0,1,\eqno{(1.10)}$$
then for $\delta'$, $\lambda$, $\varepsilon$ as above and $|\alpha_1|,|\alpha_2|\le 1$, we have the estimate 
$$\left\|\langle x\rangle^{-\frac{3+\delta'}{2}}\partial_x^{\alpha_1}\left(G-\lambda^2\pm i\varepsilon\right)^{-2}
\partial_x^{\alpha_2}
\langle x\rangle^{-\frac{3+\delta'}{2}}\right\|_{L^2\to L^2}\le C\lambda^{|\alpha_1|+|\alpha_2|-2}.\eqno{(1.11)}$$
\end{Theorem}

In fact, some of the conditions above can be weakened. 
Indeed, using Theorem 1.1 we prove the following

\begin{corol} Let $b\in L^\infty({\bf R}^n;{\bf R}^n)$, $V\in L^\infty({\bf R}^n;{\bf R})$ satisfy
$$\langle x\rangle^\delta|b(x)|+|V(x)|\le C,\,\,\forall x\in{\bf R}^n,\eqno{(1.12)}$$
with some constants $C>0$, $0<\delta\ll 1$. Suppose also that there exists a constant $r_0\gg 1$ so that $b=b^L+b^S$,
$V=V^L+V^S$ with functions $b^L,b^S\in L^\infty({\bf R}^n;{\bf R}^n)$, $V^L,V^S\in L^\infty({\bf R}^n;{\bf R})$,
$b^L$ and $V^L$ belonging to $C^1([r_0,+\infty))$ with respect to the radial variable $r$, 
and satisfying
$$\left|\partial_rb^L(rw)\right|+\left|\partial_rV^L(rw)\right|+\left|b^S(rw)\right|+\left|V^S(rw)\right|\le Cr^{-1-\delta}
\eqno{(1.13)}$$
for all $r\ge r_0$, $w\in {\bf S}^{n-1}$. Finally, we suppose that the functions $b^L(rw)$ and $b^S(rw)$ are continuous
with respect to $r$ uniformly on $[0,+\infty)\times{\bf S}^{n-1}$ and that $b(0)=0$. Then the estimate (1.7) holds true.
\end{corol}

These resolvent estimates are sharp in $\lambda$ in the sense that we have the same for the free Laplacian. 
The estimate (1.7) is well known
to hold for non-trapping compactly supported perturbations of the Laplacian 
(in which case it can be derived from the propagation of the singularities, 
e.g. see \cite{kn:V1}) and in particular when $b,V\in C_0^\infty({\bf R}^n)$, $n\ge 2$. 
It is also proved in many situations for operators of the form $-\Delta_g+V$ under the non-trapping condition, 
where $\Delta_g$ denotes the (negative) Laplace-Beltrami operator on an infinite volume unbounded Riemannian manifold 
(e.g. see \cite{kn:RT}, \cite{kn:V1}). Note that without the non-trapping condition we have in general resolvent estimates
with $O\left(e^{\gamma\lambda}\right)$, $\gamma>0$, in the right-hand side (see \cite{kn:CV1}). 
The estimate (1.7) is well known for operators $-\Delta+V$ on ${\bf R}^n$ for short-range
potentials $V\in L^\infty({\bf R}^n)$. In the case when the magnetic potential is not identically zero, 
it can also be easily proved for small short-range magnetic potentials (e.g. see \cite{kn:DF2}). 
For large short-range magnetic potentials 
$b(x)$ and electric potentials $V(x)$ the estimate (1.7) is proved in \cite{kn:EGS} (see Proposition 4.3) 
in all dimensions $n\ge 3$, 
provided $b(x)$ is a continuous function. For large long-range magnetic and electric potentials the estimate (1.7) is proved in
\cite{kn:R}, provided $b,V\in C^\infty({\bf R}^n)$ and $\partial_x^\alpha b(x), \partial_x^\alpha V(x)=O_\alpha\left(
\langle x\rangle^{-\delta-|\alpha|}\right)$, $\delta>0$. In fact, the method of \cite{kn:R} requires this condition for
$|\alpha|\le 2$, only. Note also that resolvent 
estimates like (1.7) play crucial role in the proof of uniform local energy, smoothing, Strichartz and 
dispersive estimates for the wave and the Schr\"odinger equations, which in turn explains the big interest 
in proving such kind of estimates in various situations. Therefore, the sharpness in $\lambda$ is important as a 
loss in $\lambda\gg 1$ in the resolvent estimate produces a loss of
derivatives in the applications mentioned above. 

The price to pay for assuming regularity of the potentials only with respect to the radial variable $r$ in the conditions of 
Theorem 1.1 and Corollary 1.2 is that we must also assume that the magnetic potential $b(x)$ vanishes at the origin $x=0$. In fact, we can
remove this latter condition if the regularity is assumed with respect to the variable $x$. More precisely, 
we have the following

\begin{corol}Let $b=b^L+b^S$, $V=V^L+V^S$, where $b^L\in C^1({\bf R}^n;{\bf R}^n)$, $b^S\in C^0({\bf R}^n;{\bf R}^n)$,
$V^L\in C^1({\bf R}^n;{\bf R})$, $V^S\in L^\infty({\bf R}^n;{\bf R})$ satisfy
$$\left|V^L(x)\right|+\langle x\rangle^{1+\delta}\sum_{|\alpha|=1}\left|\partial_x^\alpha V^L(x)\right|\le C,\eqno{(1.14)}$$
$$\left|V^S(x)\right|\le C\langle x\rangle^{-1-\delta},\eqno{(1.15)}$$
 $$\sum_{|\alpha|\le 1}\langle x\rangle^{|\alpha|+\delta}\left|\partial_x^\alpha b^L(x)\right|\le C,\eqno{(1.16)}$$
 $$\left|b^S(x)\right|\le C\langle x\rangle^{-1-\delta},\eqno{(1.17)}$$
$$\forall\epsilon>0\,\exists\theta=\theta(\epsilon)>0 \,\,\mbox{so that} \,\,\left|b^S(x+\theta y)-b^S(x)\right|
\le \epsilon\langle x\rangle^{-1-\delta}$$ $$\mbox{for all}\,\, x,y\in{\bf R}^{n},\,|y|\le 1,\eqno{(1.18)}$$
with some constants $C>0$ and $0<\delta\ll 1$. Then the estimate (1.7) holds true.
\end{corol}

As mentioned above, this result is proved in \cite{kn:EGS} in the case $b^L\equiv V^L\equiv 0$ by a different method.
Here we extend it to more general perturbations and provide a simpler proof.

We will use Theorem 1.1 to prove dispersive estimates for the wave group $e^{it\sqrt{G}}$ for self-adjoint operators 
$G$ as above in the case $n=3$. More precisely, we are interested in generalizing the following three dimensional dispersive estimate
$$\left\|e^{it\sqrt{G_0}}G_0^{-1-\epsilon}\chi_a(\sqrt{G_0})\right\|_{L^1\to L^\infty}\le C_{a,\epsilon}|t|^{-1},
\quad\forall t\neq 0,\eqno{(1.19)}$$
for every $a,\epsilon>0$, where $G_0$ denotes the self-adjoint realization of the free Laplacian $-\Delta$ on 
$L^2({\bf R}^3)$ and $\chi_a\in C^\infty({\bf R})$, $\chi_a(\lambda)=0$ for $\lambda\le a$, $\chi_a(\lambda)=1$ for  
$\lambda\ge a+1$. We suppose that the magnetic potential $b$ is $C^1({\bf R}^+)$ with respect to the radial variable $r$, while no 
regularity is assumed on the electric potential $V$. 
We also suppose that there exist constants $C>0$ and $0<\delta\ll 1$ such that
$$\left|V(rw)\right|+|b(rw)|\le C\langle r\rangle^{-2-\delta},\eqno{(1.20)}$$
$$\left|b(rw)\right|\le Cr^\delta\quad{\rm for}\quad r\le 1,\eqno{(1.21)}$$
$$\left|\partial_rb(rw)\right|\le Cr^{-1+\delta}\langle r\rangle^{-1-2\delta}.\eqno{(1.22)}$$
Clearly, the conditions of Theorem 1.1 are fulfilled (with $b^S\equiv V^L\equiv 0$) for $b$ and $V$ satisfying (1.20), 
(1.21) and (1.22), so the estimates (1.7) and (1.11) are valid. When $n=3$ we have the following

\begin{Theorem}Under the assumptions (1.20), (1.21) and (1.22), there exists a constant $a>0$ so that the following 
dispersive estimate holds
$$\left\|e^{it\sqrt{G}}G^{-3/2-\epsilon}\chi_a(\sqrt{G})\right\|_{L^1\to L^\infty}\le C_{\epsilon}|t|^{-1},
\quad\forall t\neq 0,\eqno{(1.23)}$$
for every $\epsilon>0$. Moreover, for every $\delta'>0$ there exists a constant $a>0$ so that we have the estimate
 $$\left\|e^{it\sqrt{G}}G^{-1-\epsilon}\chi_a(\sqrt{G})
\langle x\rangle^{-3/2-\delta'}\right\|_{L^2\to L^\infty}\le C_{\epsilon,\delta'}|t|^{-1},
\quad\forall t\neq 0,\eqno{(1.24)}$$
for every $\epsilon>0$.
\end{Theorem}

\noindent 
{\bf Remark.} In fact, one can show that the estimates (1.23) and (1.24) hold true for every $a>0$. Indeed, according
to the results of \cite{kn:KT} the condition (1.20) guarantees that the operator $G$ has no embedded strictly positive
eigenvalues, which in turn implies that the resolvent estimates (1.7) and (1.11) are valid for every $\lambda_0>0$ with
constants $C>0$ depending on $\lambda_0$.

The estimates (1.23) and (1.24) are not optimal--for example, in (1.21) there is a loss of one derivative. The desired result
would be to prove the dispersive estimate
$$\left\|e^{it\sqrt{G}}G^{-1-\epsilon}\chi_a(\sqrt{G})\right\|_{L^1\to L^\infty}\le C_{\epsilon}|t|^{-1},
\quad\forall t\neq 0,\eqno{(1.25)}$$
for every $\epsilon>0$ and some $a>0$. When $b\equiv 0$ and for a large class of rough potentials $V$ the estimate 
(1.25) follows from \cite{kn:DP}. In higher dimensions $n\ge 4$ an analogue of (1.25) is proved in 
\cite{kn:B} for Schwartz class potentials $V$ and in \cite{kn:CV2} for potentials 
$V\in C^{\frac{n-3}{2}}({\bf R}^n)$, $4\le n\le 7$, while in \cite{kn:V2} dispersive estimates with a loss of
$\frac{n-3}{2}$ derivatives are proved for potentials $V\in L^\infty({\bf R}^n)$, 
$V(x)=O\left(\langle x\rangle^{-\frac{n+1}{2}-\delta}\right)$, $\delta>0$. Proving (1.25) when the magnetic potential 
$b(x)$ is not identically zero, however, is a difficult and an open problem even if $b$ is supposed small and smooth. 
Our conjecture is that (1.25) should hold for $b\in C_0^1({\bf R}^3)$
and $V\in L^\infty({\bf R}^3)$, $V(x)=O\left(\langle x\rangle^{-2-\delta}\right)$, $\delta>0$, 
while in higher dimensions $n\ge 4$ we expect 
to have an optimal dispersive estimate (that is, without loss of derivatives) similar to (1.25) for 
$b\in C_0^{\frac{n-1}{2}}({\bf R}^n)$
and $V\in C_0^{\frac{n-3}{2}}({\bf R}^n)$. Note that dispersive estimates for the wave group with a loss of 
$\frac{n}{2}$ derivatives have been
recently proved in \cite{kn:CCV} in all dimensions $n\ge 2$ for a class of potentials $b\in C^1({\bf R}^n)$ and 
$V\in L^\infty({\bf R}^n)$. 
Note also that an estimate similar to (1.24) is proved in \cite{kn:DF} for a class of small
potentials $b$ and $V$ still in dimension three. 

Theorem 1.1 plays a crucial role in the proof of the dispersive estimates (1.23) and (1.24). Note that we cannot use 
Corollary 1.3 instead, since a function $b(x)$ satisfying the conditions (1.20), (1.21) and
(1.22) is not necessarily continuous in $x$. Finally, we expect that Theorem 1.4 can be extended to all dimensions
$n\ge 3$ for potentials $b(x), V(x)=O\left(\langle x\rangle^{-\frac{n+1}{2}-\delta}\right)$.

\section{Resolvent estimates}
Clearly, it suffices to prove the resolvent estimates for $0<\delta'\le\delta$. We will first consider the case 
$b^S\equiv V^S\equiv 0$, so $b=b^L$ and $V=V^L$. Let $\alpha_1=\alpha_2=0$. Clearly, in this case (1.7) follows 
from the a priori estimate
$$\left\|\psi_{\delta'}(|x|)^{1/2}f\right\|_{L^2({\bf R}^n)}\le C\lambda^{-1}
\left\|\psi_{\delta'}(|x|)^{-1/2}\left(
G-\lambda^2\pm i\varepsilon\right)f\right\|_{L^2({\bf R}^n)}.\eqno{(2.1)}$$
It suffices to consider the case ``+`` only. 
To prove (2.1) we will pass to polar coordinates $(r,w)\in {\bf R}^+\times{\bf S}^{n-1}$.
Recall that $L^2({\bf R}^n)\cong L^2\left({\bf R}^+\times{\bf S}^{n-1},r^{(n-1)/2}drdw\right)$. Set $X=
\left({\bf R}^+\times{\bf S}^{n-1},drdw\right)$, $u=r^{(n-1)/2}f$, 
$$P=\lambda^{-2}r^{(n-1)/2}\left(G-\lambda^2+i\varepsilon\right)r^{-(n-1)/2}.$$
It is well known that
$$r^{(n-1)/2}\Delta r^{-(n-1)/2}=\partial_r^2+\frac{\Delta_w-c_n}{r^2},\eqno{(2.2)}$$
where 
$$c_n=\frac{(n-1)(n-3)}{4}$$
and $\Delta_w$ denotes the (negative) Laplace-Beltrami operator on ${\bf S}^{n-1}$ written in the coordinates $w$.
It is easy to see that (2.1) follows from the estimate
$$\left\|\psi_{\delta'}(r)^{1/2}u\right\|_{H^1(X)}\le C\lambda\left\|\psi_{\delta'}(r)^{-1/2}Pu\right\|_{L^2(X)},
\eqno{(2.3)}$$
where the norm in the left-hand side is defined as follows
$$\left\|\psi_{\delta'}(r)^{1/2}u\right\|_{H^1(X)}^2=\left\|\psi_{\delta'}(r)^{1/2}u\right\|_{L^2(X)}^2+
\left\|\psi_{\delta'}(r)^{1/2}{\cal D}_ru\right\|_{L^2(X)}^2+\left\|\psi_{\delta'}(r)^{1/2}r^{-1}
\Lambda_w^{1/2}u\right\|_{L^2(X)}^2,$$
where ${\cal D}_r=i\lambda^{-1}\partial_r$, $\Lambda_w=-\lambda^{-2}\Delta_w$. 
Througout this section $\|\cdot\|$ and $\langle\cdot,\cdot\rangle$ will denote the norm and the scalar product in the
Hilbert space $L^2({\bf S}^{n-1})$. Hence $\|u\|^2_{L^2(X)}=\int_0^\infty\|u(r,\cdot)\|^2dr$.
 Using (2.2) one can easily check that the operator $P$ can be written in the form
$$P={\cal D}_r^2+r^{-2}\widetilde\Lambda_w+\lambda^{-2}W(r,w)-1+i\varepsilon\lambda^{-2}$$ $$+\lambda^{-1}\sum_{j=1}^nw_j
\left(b_{j}(rw){\cal D}_r+{\cal D}_rb_{j}(rw)\right)$$ $$+\lambda^{-1}r^{-1}
\sum_{j=1}^n\left(b_{j}(rw)Q_j(w,{\cal D}_w)+
Q_j(w,{\cal D}_w)b_{j}(rw)\right),$$
 $$W=V(rw)+|b(rw)|^2-i(n-1)r^{-1}\sum_{j=1}^nw_jb_j(rw),$$
where $\widetilde\Lambda_w=\Lambda_w+\lambda^{-2}c_n$, $w_j=x_j/r$, ${\cal D}_w=i\lambda^{-1}\partial_w$, 
$Q_j(w,{\cal D}_w)=i\lambda^{-1}Q_j(w,\partial_w)$, $Q_j(w,\xi)\in C^\infty(T^*{\bf S}^{n-1})$ 
are real-valued, independent of $r$ and $\lambda$, and homogeneous of order 1 with respect to $\xi$.
Decompose $W$ as $W^L+W^S$, where
$$W^L=V(rw)+\left|b(rw)\right|^2,$$
$$W^S=-i(n-1)r^{-1}\sum_{j=1}^nw_jb_j(rw).$$
It is easy to see that the assumptions (1.1), (1.2) and (1.4) imply
$$\left|W^L(r,w)\right|\le C,\eqno{(2.4)}$$
$$\partial_rW^L(r,w)\le C\psi_\delta(r),\eqno{(2.5)}$$
$$\left|W^S(r,w)\right|\le C\psi_\delta(r).\eqno{(2.6)}$$
Set 
$$E(r)=-\left\langle\left(r^{-2}\widetilde\Lambda_w-1+\lambda^{-2}W^L\right)u(r,w),u(r,w)\right\rangle+
\left\|{\cal D}_ru(r,w)\right\|^2$$
  $$-2\lambda^{-1}r^{-1}\sum_{j=1}^n{\rm Re}\,\left\langle b_{j}(rw)Q_j(w,{\cal D}_w)u(r,w),u(r,w)\right\rangle.$$
 We have the identity
$$E'(r):=\frac{dE(r)}{dr}=\frac{2}{r}\left\langle r^{-2}\widetilde\Lambda_wu(r,w),u(r,w)
\right\rangle-\lambda^{-2}\left\langle 
\frac{\partial W^L}{\partial r}u(r,w),u(r,w)\right\rangle$$
 $$-2\lambda^{-1}\sum_{j=1}^n{\rm Re}\,\left\langle\frac{\partial(b_{j}(rw)/r)}{\partial r}
Q_j(w,{\cal D}_w)u(r,w),u(r,w)\right\rangle$$
 $$-2\lambda^{-1}\sum_{j=1}^n{\rm Re}\,\left\langle w_j\frac{\partial b_{j}(rw)}{\partial r}
u(r,w),{\cal D}_ru(r,w)\right\rangle
+2\lambda{\rm Im}\,\left\langle \widetilde Pu(r,w),{\cal D}_ru(r,w)\right\rangle,$$
where
$$\widetilde P=P-i\varepsilon\lambda^{-2}-\lambda^{-2}W^S.$$
Observe now that by (1.4) we have
$$\left|\frac{\partial b(rw)}{\partial r}\right|\le C\psi_\delta(r),\eqno{(2.7)}$$
$$\left|\frac{\partial(b(rw)/r)}{\partial r}\right|\le Cr^{-3/2}\psi_\delta(r)^{1/2}.\eqno{(2.8)}$$
Hence, using (2.5), (2.7) and (2.8), we obtain
$$E'(r)\ge\frac{2}{r}\left\langle r^{-2}\widetilde\Lambda_wu(r,w),u(r,w)\right\rangle-\gamma r^{-3}\sum_{j=1}^n\left\|
Q_j(w,{\cal D}_w)u(r,w)\right\|^2$$ $$-\lambda^{-1}\left\|\psi_\delta^{1/2}{\cal D}_ru(r,w)\right\|^2-O_{\gamma}
(\lambda^{-1})\left\|\psi_\delta^{1/2}u(r,w)\right\|^2-2\lambda M(r),\eqno{(2.9)}$$
$\forall \gamma>0$ independent of $\lambda$ and $r$, where
$$M(r)= \left|\left\langle \widetilde Pu(r,w),{\cal D}_ru(r,w)\right\rangle\right|.$$
Since $\|Q_j(w,{\cal D}_w)u\|\le C\|\Lambda_w^{1/2}u\|\le C\|\widetilde\Lambda_w^{1/2}u\|$, taking $\gamma$ small enough we can absorb the second term in the 
right-hand side of (2.9) by the first one and obtain
$$E'(r)\ge\frac{1}{r}\left\langle r^{-2}\widetilde\Lambda_wu(r,w),u(r,w)\right\rangle-\lambda^{-1}
\left\|\psi_\delta^{1/2}{\cal D}_ru(r,w)\right\|^2$$ $$-
O(\lambda^{-1})\left\|\psi_\delta^{1/2}u(r,w)\right\|^2-2\lambda M(r).\eqno{(2.10)}$$
Using that $\widetilde\Lambda_w\ge 0$, we deduce from (2.10)
$$E(r)=-\int_r^\infty E'(t)dt\le \lambda^{-1}\left\|\psi_\delta^{1/2}{\cal D}_ru\right\|^2_{L^2(X)}+
O(\lambda^{-1})\left\|\psi_\delta^{1/2}u\right\|^2_{L^2(X)}+2\lambda \int_0^\infty M(t)dt.\eqno{(2.11)}$$
Let now $\psi(r)>0$ be such that $\int_0^\infty\psi(r)dr<+\infty$. Multiplying both sides of (2.11) by $\psi$ and integrating
from $0$ to $\infty$, we get
$$\int_0^\infty\psi(r)E(r)dr\le O(\lambda^{-1})\left\|\psi_\delta^{1/2}{\cal D}_ru\right\|^2_{L^2(X)}+
O(\lambda^{-1})\left\|\psi_\delta^{1/2}u\right\|^2_{L^2(X)}+O(\lambda)\int_0^\infty M(r)dr.\eqno{(2.12)}$$
In particular, (2.12) holds with $\psi=\psi_{\delta'}(r)$ for any $0<\delta'\le\delta$. 
It is easy also to check that
$$0<-\frac{d}{dr}\left(r\psi_{\delta'}(r)\right)\le C\psi_{\delta'}(r),$$
so we can use (2.12) with $\psi=-\frac{d}{dr}\left(r\psi_{\delta'}(r)\right)$ to obtain
$$\int_0^\infty r\psi_{\delta'}(r)E'(r)dr=-\int_0^\infty\frac{d}{dr}\left(r\psi_{\delta'}(r)\right)E(r)dr$$
$$\le O(\lambda^{-1})\left\|\psi_\delta^{1/2}{\cal D}_ru\right\|^2_{L^2(X)}+
O(\lambda^{-1})\left\|\psi_\delta^{1/2}u\right\|^2_{L^2(X)}+O(\lambda)\int_0^\infty M(r)dr.\eqno{(2.13)}$$
Since $r\psi_{\delta'}(r)\le 1$, combining (2.10) and (2.13) we conclude
$$\left\|\psi_{\delta'}^{1/2}r^{-1}\widetilde\Lambda_w^{1/2}u\right\|^2_{L^2(X)}
\le O(\lambda^{-1})\left\|\psi_\delta^{1/2}{\cal D}_ru\right\|^2_{L^2(X)}$$ $$+
O(\lambda^{-1})\left\|\psi_\delta^{1/2}u\right\|^2_{L^2(X)}+O(\lambda)\int_0^\infty M(r)dr.\eqno{(2.14)}$$
On the other hand, in view of (2.4) we can choose $\lambda$ big enough so that $1-\lambda^{-2}W^L\ge 1/2$.
Therefore, for $\lambda\gg 1$ we have the inequality
$$\int_0^\infty\psi_{\delta'}(r)E(r)dr\ge \left\|\psi_{\delta'}^{1/2}{\cal D}_ru\right\|^2_{L^2(X)}
+\frac{1}{3}\left\|\psi_{\delta'}^{1/2}u\right\|^2_{L^2(X)}-2\left\|\psi_{\delta'}^{1/2}r^{-1}\widetilde\Lambda_w^{1/2}u
\right\|^2_{L^2(X)}.\eqno{(2.15)}$$
By (2.12), (2.14) and (2.15), we conclude
$$\left\|\psi_{\delta'}^{1/2}u\right\|^2_{\widetilde H^1(X)}:=\left\|\psi_{\delta'}^{1/2}{\cal D}_ru\right\|^2_{L^2(X)}
+\left\|\psi_{\delta'}^{1/2}u\right\|^2_{L^2(X)}+\left\|\psi_{\delta'}^{1/2}r^{-1}\widetilde\Lambda_w^{1/2}u
\right\|^2_{L^2(X)}$$
 $$\le O(\lambda^{-1})\left\|\psi_\delta^{1/2}{\cal D}_ru\right\|^2_{L^2(X)}+
O(\lambda^{-1})\left\|\psi_\delta^{1/2}u\right\|^2_{L^2(X)}+O(\lambda) \int_0^\infty M(r)dr.\eqno{(2.16)}$$
Set
$$P^\sharp=\widetilde P+i\varepsilon\lambda^{-2}=P-\lambda^{-2}W^S(r,w),$$
 $$M^\sharp(r)= \left|\left\langle P^\sharp u(r,w),{\cal D}_ru(r,w)\right\rangle\right|,
\quad N(r)= \left|\left\langle Pu(r,w),{\cal D}_ru(r,w)\right\rangle\right|.$$
In view of (2.6), we have
$$\lambda \int_0^\infty M^\sharp(r)dr\le \lambda\int_0^\infty N(r)dr+O(\lambda^{-1})
\left\|\psi_\delta^{1/2}u\right\|^2_{L^2(X)}.\eqno{(2.17)}$$
We also have
$$\lambda \int_0^\infty M(r)dr\le \lambda \int_0^\infty M^\sharp(r)dr+\varepsilon\lambda^{-1}\left(\|u\|_{L^2(X)}^2
+\left\|{\cal D}_ru\right\|_{L^2(X)}^2\right),\eqno{(2.18)}$$ 
$$\lambda\int_0^\infty N(r)dr\le O_\gamma(\lambda^2)\left\|\psi_{\delta'}^{-1/2}Pu\right\|_{L^2(X)}^2+\gamma
\left\|\psi_{\delta'}^{1/2}{\cal D}_ru\right\|_{L^2(X)}^2,\eqno{(2.19)}$$
for every $\gamma>0$ independent of $\lambda$. 
On the other hand, in view of (2.4) and (2.6), we have
$$\varepsilon\lambda^{-2}\|u\|_{L^2(X)}^2={\rm Im}\,\left\langle Pu,u\right\rangle_{L^2(X)}+(n-1)\lambda^{-2}\sum_{j=1}^n
\left\langle r^{-1}w_jb_j(rw)u,u\right\rangle_{L^2(X)}$$ $$
\le \left|\left\langle Pu,u\right\rangle_{L^2(X)}\right|+O(\lambda^{-2})\left\|\psi_\delta^{1/2}
u\right\|_{L^2(X)}^2,\eqno{(2.20)}$$
 $${\rm Re}\,\left\langle Pu,u\right\rangle_{L^2(X)}=\left\|{\cal D}_ru\right\|_{L^2(X)}^2+
\left\|r^{-1}\widetilde\Lambda_w^{1/2}u\right\|_{L^2(X)}^2+\left\langle\left(\lambda^{-2}W^L-1\right)u,u
\right\rangle_{L^2(X)}$$
 $$+2\lambda^{-1}\sum_{j=1}^n{\rm Re}\,\left\langle w_jb_{j}(rw){\cal D}_ru,u\right\rangle_{L^2(X)}+
2\lambda^{-1}\sum_{j=1}^n{\rm Re}\,\left\langle r^{-1}b_{j}(rw)Q_ju,u\right\rangle_{L^2(X)}$$ 
 $$\ge \left\|{\cal D}_ru\right\|_{L^2(X)}^2+
\left\|r^{-1}\widetilde\Lambda_w^{1/2}u\right\|_{L^2(X)}^2-O(1)\left\|u\right\|_{L^2(X)}^2$$ $$-O(\lambda^{-1})\left(
\left\|{\cal D}_ru\right\|_{L^2(X)}^2+
\left\|r^{-1}\widetilde\Lambda_w^{1/2}u\right\|_{L^2(X)}^2+\left\|u\right\|_{L^2(X)}^2\right)$$
 $$\ge \frac{1}{2}\left\|{\cal D}_ru\right\|_{L^2(X)}^2-O(1)\left\|u\right\|_{L^2(X)}^2,$$
provided $\lambda$ is taken large enough, which in turn implies
$$\left\|{\cal D}_ru\right\|_{L^2(X)}^2\le O(1)\left\|u\right\|_{L^2(X)}^2+
2\left|\left\langle Pu,u\right\rangle_{L^2(X)}\right|.\eqno{(2.21)}$$
Combining (2.18), (2.20) and (2.21), we get
$$\lambda \int_0^\infty M(r)dr\le O(\lambda)\int_0^\infty M^\sharp(r)dr+O(\lambda^{-1})\left\|\psi_\delta^{1/2}
u\right\|_{L^2(X)}^2+O(\lambda)
\left|\left\langle Pu,u\right\rangle_{L^2(X)}\right|$$ $$\le O(\lambda)\int_0^\infty M^\sharp(r)dr+
O_\gamma(\lambda^2)\left\|\psi_{\delta'}^{-1/2}Pu\right\|_{L^2(X)}^2+O(\gamma+\lambda^{-1})
\left\|\psi_{\delta'}^{1/2}u\right\|_{L^2(X)}^2.\eqno{(2.22)}$$
By (2.16), (2.17), (2.19) and (2.22), we conclude
$$\left\|\psi_{\delta'}^{1/2}u\right\|_{\widetilde H^1(X)}\le \left(O_\gamma(\lambda^{-1})
+O(\gamma)\right)^{1/2}\left\|\psi_{\delta'}^{1/2}u
\right\|_{\widetilde H^1(X)}+O_\gamma(\lambda)\left\|\psi_{\delta'}^{-1/2}Pu\right\|_{L^2(X)},\eqno{(2.23)}$$
where we have used that $\psi_\delta\le\psi_{\delta'}$ for $\delta'\le\delta$. Now, taking $\gamma$ small enough,
independent of $\lambda$, and $\lambda$ big enough, we can absorb the first term in the right-hand side of (2.23) to
obtain (2.3). To prove (1.7) for all multi-indices $|\alpha_1|,|\alpha_2|\le 1$ we will use the following

\begin{lemma} If $|b(x)|+|V(x)|\le C'=Const$, then for every $s\in{\bf R}$ there exist constants $C>0$ independent of 
$b$ and $V$ and $\lambda_0>0$ depending on $C'$ so that for $\lambda\ge\lambda_0$ and $0\le|\alpha_1|,|\alpha_2|\le 1$ 
we have the estimate
 $$\left\|\langle x\rangle^{-s}\partial_x^{\alpha_1}\left(G\pm i\lambda^2\right)^{-1}\partial_x^{\alpha_2}
\langle x\rangle^{s}\right\|_{L^2\to L^2}\le C\lambda^{|\alpha_1|+|\alpha_2|-2}.\eqno{(2.24)}$$
\end{lemma}

{\it Proof.} Without loss of generality we may suppose that $s\ge 0$. Let us first see that (2.24) is valid for the free
operator $G_0$. This is obvious for $s=0$. For $s>0$ we will use the identity
$$\left(G_0\pm i\lambda^2\right)^{-1}\langle x\rangle^{s}=\langle x\rangle^{s}\left(G_0\pm i\lambda^2\right)^{-1}-
\left(G_0\pm i\lambda^2\right)^{-1}\left[\Delta,\langle x\rangle^{s}\right]\left(G_0\pm i\lambda^2\right)^{-1}.\eqno{(2.25)}$$
Since
$$\left[\Delta,\langle x\rangle^{s}\right]=O\left(\langle x\rangle^{s-1}\right)\partial_x+
O\left(\langle x\rangle^{s-2}\right),$$
we obtain from (2.25) (with $|\alpha|\le 2$)
$$\left\|\langle x\rangle^{-s}\partial_x^\alpha \left(G_0\pm i\lambda^2\right)^{-1}\langle x\rangle^{s}
\right\|_{L^2\to L^2}\le \left\|\langle x\rangle^{-s}\partial_x^\alpha \langle x\rangle^{s}\left(G_0\pm i\lambda^2\right)^{-1}
\right\|_{L^2\to L^2}$$
 $$+C\sum_{|\beta|\le 1}\left\|\langle x\rangle^{-s}\partial_x^\alpha \left(G_0\pm i\lambda^2\right)^{-1}\langle x\rangle^{s-1}
\right\|_{L^2\to L^2}\left\|\partial_x^\beta\left(G_0\pm i\lambda^2\right)^{-1}\right\|_{L^2\to L^2}$$
 $$\le C\lambda^{|\alpha|-2}+O(\lambda^{-1})\left\|\langle x\rangle^{-s}\partial_x^\alpha \left(G_0\pm i\lambda^2\right)^{-1}\langle x\rangle^{s-1}
\right\|_{L^2\to L^2}.\eqno{(2.26)}$$
Iterating (2.26) a finite number of times and taking into account that the operator $\partial_x^{\alpha_2}$ commutes
with the free resolvent, we get (2.24) for $G_0$. To prove (2.24) for the perturbed operator we will use the resolvent
identity
$$\left(G\pm i\lambda^2\right)^{-1}=\left(G_0\pm i\lambda^2\right)^{-1}-\left(G\pm i\lambda^2\right)^{-1}(G-G_0)
\left(G_0\pm i\lambda^2\right)^{-1}.\eqno{(2.27)}$$
By (2.27) we get
$$\sum_{|\alpha_1|,|\alpha_2|\le 1}\lambda^{-|\alpha_1|-|\alpha_2|}\left\|\langle x\rangle^{-s}\partial_x^{\alpha_1}
\left(G\pm i\lambda^2\right)^{-1}\partial_x^{\alpha_2}\langle x\rangle^{s}\right\|_{L^2\to L^2}$$ 
 $$\le\sum_{|\alpha_1|,|\alpha_2|\le 1}\lambda^{-|\alpha_1|-|\alpha_2|}\left\|\langle x\rangle^{-s}\partial_x^{\alpha_1}
\left(G_0\pm i\lambda^2\right)^{-1}\partial_x^{\alpha_2}\langle x\rangle^{s}\right\|_{L^2\to L^2}$$ 
 $$+C\sum_{|\alpha_1|,|\alpha_2|\le 1}\sum_{|\beta_1|+|\beta_2|\le 1}\lambda^{-|\alpha_1|-|\alpha_2|}
\left\|\langle x\rangle^{-s}\partial_x^{\alpha_1}
\left(G\pm i\lambda^2\right)^{-1}\partial_x^{\beta_1}\langle x\rangle^{s}\right\|_{L^2\to L^2}$$ 
$$\times \left\|\langle x\rangle^{-s}\partial_x^{\beta_2}
\left(G_0\pm i\lambda^2\right)^{-1}\partial_x^{\alpha_2}\langle x\rangle^{s}\right\|_{L^2\to L^2}$$
 $$\le C\lambda^{-2}+O(\lambda^{-1})\sum_{|\alpha_1|,|\beta_1|\le 1}\lambda^{-|\alpha_1|-|\beta_1|}
\left\|\langle x\rangle^{-s}\partial_x^{\alpha_1}
\left(G\pm i\lambda^2\right)^{-1}\partial_x^{\beta_1}\langle x\rangle^{s}\right\|_{L^2\to L^2}.\eqno{(2.28)}$$ 
Taking now $\lambda$ big enough we can absorb the second term in the right-hand side of (2.28) and obtain (2.24).
\eproof

Let us see that (1.7) for all multi-indices $\alpha_1$ and $\alpha_2$ follows from (1.7) 
with $\alpha_1=\alpha_2=0$ and Lemma 2.1. To this end, we will use the resolvent identity
$$\left(G-\lambda^2\pm i\varepsilon\right)^{-1}=\left(G-i\lambda^2\right)^{-1}+(\lambda^2\mp i\varepsilon-i\lambda^2)
\left(G-i\lambda^2\right)^{-2}$$ $$+(\lambda^2\mp i\varepsilon-i\lambda^2)^2\left(G-i\lambda^2\right)^{-1}
\left(G-\lambda^2\pm i\varepsilon\right)^{-1}\left(G-i\lambda^2\right)^{-1}.$$
Hence
$$\left\|\langle x\rangle^{-\frac{1+\delta'}{2}}\partial_x^{\alpha_1}\left(G-\lambda^2\pm i\varepsilon\right)^{-1}
\partial_x^{\alpha_2}\langle x\rangle^{-\frac{1+\delta'}{2}}\right\|_{L^2\to L^2}\le \left\|\partial_x^{\alpha_1}
\left(G-i\lambda^2\right)^{-1}\partial_x^{\alpha_2}\right\|_{L^2\to L^2}$$
 $$+C\lambda^2\left\|\partial_x^{\alpha_1}\left(G-i\lambda^2\right)^{-1}\right\|_{L^2\to L^2}\left\|
\left(G-i\lambda^2\right)^{-1}\partial_x^{\alpha_2}\right\|_{L^2\to L^2}$$
 $$+C\lambda^4\left\|\langle x\rangle^{-\frac{1+\delta'}{2}}\partial_x^{\alpha_1}\left(G-i\lambda^2\right)^{-1}
\langle x\rangle^{\frac{1+\delta'}{2}}\right\|_{L^2\to L^2}$$ $$\times
\left\|\langle x\rangle^{-\frac{1+\delta'}{2}}\left(G-\lambda^2\pm i\varepsilon\right)^{-1}
\langle x\rangle^{-\frac{1+\delta'}{2}}\right\|_{L^2\to L^2}$$ 
$$\times \left\|\langle x\rangle^{\frac{1+\delta'}{2}}\left(G-i\lambda^2\right)^{-1}\partial_x^{\alpha_2}
\langle x\rangle^{-\frac{1+\delta'}{2}}\right\|_{L^2\to L^2}$$
 $$\le C\lambda^{|\alpha_1|+|\alpha_2|-2}+C\lambda^{|\alpha_1|+|\alpha_2|}\left\|\langle x\rangle^{-\frac{1+\delta'}{2}}
\left(G-\lambda^2\pm i\varepsilon\right)^{-1}
\langle x\rangle^{-\frac{1+\delta'}{2}}\right\|_{L^2\to L^2}$$ $$\le \widetilde C\lambda^{|\alpha_1|+|\alpha_2|-1}.$$ 
We will now prove (1.7) in the general case. 
Let $\phi\in C_0^\infty({\bf R}^+)$, $\phi\ge 0$, $\int\phi(\sigma)d\sigma=1$, and given any $0<\theta\le 1$, set
$$B_\theta(r,w)=\theta^{-1}\eta_\delta(r)\int_{{\bf R}} g_\delta(r',w)\phi\left(\frac{r'-r}{\theta}\right)dr'=
\eta_\delta(r)\int_{{\bf R}} g_\delta(r+\theta \sigma,w)\phi(\sigma)d\sigma,$$
$b_\theta^S(x):=B_\theta(|x|,\frac{x}{|x|})$. 
In view of the assumption (1.6), given any $\epsilon>0$ there exists $\theta>0$ so that for all $x\in {\bf R}^n$ we have
$$\left|b_\theta^S(x)-b^S(x)\right|\le\eta_\delta(|x|)\int_{{\bf R}}\left|g_\delta(|x|+\theta\sigma,\frac{x}{|x|})-
g_\delta(|x|,\frac{x}{|x|})\right|
\phi(\sigma)d\sigma\le \epsilon\eta_\delta(|x|).\eqno{(2.29)}$$
It is also clear that (1.5) implies the bounds
$$\left|b_\theta^S(rw)\right|\le C\eta_\delta(r),\eqno{(2.30)}$$
$$\left|\partial_rb_\theta^S(rw)\right|\le C_\epsilon \psi_\delta(r).\eqno{(2.31)}$$
We will use the above analysis and the easy observation that the constant $C$ in the right-hand side of (1.7) depends
only on the parameter $\delta'$, provided $0<\delta'\le\delta$. In view of (1.1), (1.2), (1.4), (2.30) and (2.31),
we can apply the already proved estimate (1.7) to the operator
$$G_1=-\Delta+i(b^L+b_\theta^S)\cdot\nabla +i\nabla\cdot(b^L+b_\theta^S)+V^L+\left|b^L\right|^2$$
to get the estimate
$$\left\|\langle x\rangle^{-\frac{1+\delta'}{2}}\partial_x^{\alpha_1}\left(G_1-\lambda^2\pm i\varepsilon\right)^{-1}
\partial_x^{\alpha_2}
\langle x\rangle^{-\frac{1+\delta'}{2}}\right\|_{L^2\to L^2}\le C\lambda^{|\alpha_1|+|\alpha_2|-1}\eqno{(2.32)}$$
for $\lambda\ge\lambda_0(\epsilon)>0$ with a constant $C>0$ independent of $\epsilon$, $\varepsilon$ and $\lambda$.
On the other hand, in view of (1.3), (1.5) and (2.29), the difference $G-G_1$ is a first order differential operator of
the form
$$G-G_1=O\left(\epsilon\langle x\rangle^{-1-\delta}\right)\cdot\nabla+
\nabla\cdot O\left(\epsilon\langle x\rangle^{-1-\delta}\right)+O\left(\langle x\rangle^{-1-\delta}\right).$$
Using this together with (2.32) and the resolvent identity
$$(G-\lambda^2\pm i\varepsilon)^{-1}=(G_1-\lambda^2\pm i\varepsilon)^{-1}-(G_1-\lambda^2\pm i\varepsilon)^{-1}
(G-G_1)(G-\lambda^2\pm i\varepsilon)^{-1},$$
we obtain
$$\left\|\langle x\rangle^{-\frac{1+\delta'}{2}}\left(G-\lambda^2\pm i\varepsilon\right)^{-1}
\langle x\rangle^{-\frac{1+\delta'}{2}}\right\|_{L^2\to L^2}$$ $$\le \left\|\langle x\rangle^{-\frac{1+\delta'}{2}}
\left(G_1-\lambda^2\pm i\varepsilon\right)^{-1}\langle x\rangle^{-\frac{1+\delta'}{2}}\right\|_{L^2\to L^2}$$
 $$+C\sum_{|\beta_1|+|\beta_2|\le 1}\epsilon^{|\beta_1|+|\beta_2|}\left\|\langle x\rangle^{-\frac{1+\delta'}{2}}
\left(G_1-\lambda^2\pm i\varepsilon\right)^{-1}\partial_x^{\beta_1}\langle x\rangle^{-\frac{1+\delta}{2}}\right\|_{L^2\to L^2}$$
 $$\times\left\|\langle x\rangle^{-\frac{1+\delta}{2}}\partial_x^{\beta_2}\left(G-\lambda^2\pm i\varepsilon\right)^{-1}
\langle x\rangle^{-\frac{1+\delta'}{2}}\right\|_{L^2\to L^2}$$
 $$\le C\lambda^{-1}+C\sum_{|\beta_1|+|\beta_2|\le 1}\epsilon^{|\beta_1|+|\beta_2|}\lambda^{|\beta_1|-1}
\left\|\langle x\rangle^{-\frac{1+\delta}{2}}\partial_x^{\beta_2}\left(G-\lambda^2\pm i\varepsilon\right)^{-1}
\langle x\rangle^{-\frac{1+\delta'}{2}}\right\|_{L^2\to L^2}$$
  $$\le C\lambda^{-1}+O\left(\epsilon+\lambda^{-1}\right)\left\|\langle x\rangle^{-\frac{1+\delta}{2}}
\left(G-\lambda^2\pm i\varepsilon\right)^{-1}\langle x\rangle^{-\frac{1+\delta'}{2}}\right\|_{L^2\to L^2}$$
 $$+O(\epsilon\lambda^{-1})\sum_{|\beta_2|=1}\left\|\langle x\rangle^{-\frac{1+\delta}{2}}\partial_x^{\beta_2}
(G-i\lambda^2)^{-1}
\langle x\rangle^{\frac{1+\delta}{2}}\right\|_{L^2\to L^2}$$ $$\times\left(1+\lambda^2
\left\|\langle x\rangle^{-\frac{1+\delta}{2}}
\left(G-\lambda^2\pm i\varepsilon\right)^{-1}\langle x\rangle^{-\frac{1+\delta'}{2}}\right\|_{L^2\to L^2}\right)$$
 $$\le C\lambda^{-1}+O\left(\epsilon+\lambda^{-1}\right)\left\|\langle x\rangle^{-\frac{1+\delta'}{2}}
\left(G-\lambda^2\pm i\varepsilon\right)^{-1}\langle x\rangle^{-\frac{1+\delta'}{2}}\right\|_{L^2\to L^2},\eqno{(2.33)}$$ 
where we have used that $\delta'\le\delta$ and Lemma 2.1. Taking $\epsilon>0$ small enough, independent of $\lambda$,
and $\lambda$ big enough we can absorb the second term in the right-hand side of (2.33) and obtain (1.7) in the general case
when $\alpha_1=\alpha_2=0$. For all multi-indices $\alpha_1$ and $\alpha_2$ the estimate (1.7) follows from (1.7) with
$\alpha_1=\alpha_2=0$ and Lemma 2.1 in the same way as above.

To prove (1.11) we will use the commutator identity
$$\partial_r^2+\frac{\Delta_w-c_n}{r^2}+\frac{1}{2}\left[r\partial_r,\partial_r^2+\frac{\Delta_w-c_n}{r^2}\right]=0.\eqno{(2.34)}$$
We obtain from (2.34) that the operators $\widetilde G=r^{(n-1)/2}Gr^{-(n-1)/2}$ and 
$\widetilde\Delta=r^{(n-1)/2}\Delta r^{-(n-1)/2}$ satisfy the identity
$$\widetilde G+\frac{1}{2}\left[r\partial_r,\widetilde G\right]=\widetilde G+\widetilde\Delta+
\frac{1}{2}\left[r\partial_r,\widetilde G+\widetilde\Delta\right]:={\cal Q}.\eqno{(2.35)}$$
We rewrite (2.35) as follows
$$\widetilde G-\lambda^2+i\varepsilon+\frac{1}{2}\left[r\partial_r,\widetilde G-\lambda^2+i\varepsilon\right]
=-\lambda^2+i\varepsilon+{\cal Q},$$
which yields the identity
$$\left(\widetilde G-\lambda^2+i\varepsilon\right)^{-1}-\frac{1}{2}\left[r\partial_r,\left(\widetilde G-\lambda^2+
i\varepsilon\right)^{-1}\right]$$ $$=(-\lambda^2+i\varepsilon)\left(\widetilde G-\lambda^2+i\varepsilon\right)^{-2}+
\left(\widetilde G-\lambda^2+i\varepsilon\right)^{-1}{\cal Q}\left(\widetilde G-\lambda^2+i\varepsilon\right)^{-1}.\eqno{(2.36)}$$
Set
$$\widetilde W^L=V^L(rw)+\left|b(rw)\right|^2,$$
$$\widetilde W^S=V^S(rw)-i(n-1)r^{-1}\sum_{j=1}^nw_jb_j(rw).$$
Observe now that 
$${\cal Q}=\frac{1}{2r}\frac{\partial(r^2\widetilde W^L)}{\partial r}+\frac{1}{2}\left(\widetilde W^S+
\partial_rr\widetilde W^S-r\widetilde W^S\partial_r\right)$$
 $$+\frac{i}{2}
\sum_{j=1}^nw_j\left(\frac{\partial(rb_j)}{\partial r}\partial_r+\partial_r\frac{\partial(rb_j)}{\partial r}\right)+
\frac{i}{2r}\sum_{j=1}^n\left(\frac{\partial(rb_j)}{\partial r}Q_j(w,\partial_w)+Q_j(w,\partial_w)
\frac{\partial(rb_j)}{\partial r}\right).$$
It follows from the assumptions (1.4), (1.8), (1.9) and (1.10) that
$$\left|\frac{1}{r}\frac{\partial(r^2\widetilde W^L)}{\partial r}\right|+\left|\frac{\partial(rb)}{\partial r}\right|+
\langle r\rangle\left|\widetilde W^S\right|\le C\psi_\delta(r).\eqno{(2.37)}$$
By (2.36) and (2.37) we obtain
$$\lambda^2\left\|\langle r\rangle^{-1}\psi_{\delta'}(r)^{1/2}\left(\widetilde G-\lambda^2+i\varepsilon\right)^{-2}
\psi_{\delta'}(r)^{1/2}
\langle r\rangle^{-1}\right\|_{L^2(X)\to L^2(X)}$$ $$\le \left\|\psi_{\delta'}(r)^{1/2}\left(\widetilde G-\lambda^2+
i\varepsilon\right)^{-1}\psi_{\delta'}(r)^{1/2}
\right\|_{L^2(X)\to L^2(X)}$$ $$+O(\lambda)\sum_{\pm}\left\|\psi_{\delta'}(r)^{1/2}{\cal D}_r
\left(\widetilde G-\lambda^2\pm i\varepsilon\right)^{-1}\psi_{\delta'}(r)^{1/2}
\right\|_{L^2(X)\to L^2(X)}$$ $$+O(\lambda)\sum_{\pm}\left\|\psi_{\delta'}(r)^{1/2}\left(\widetilde G-\lambda^2\pm 
i\varepsilon\right)^{-1}\psi_{\delta'}(r)^{1/2}
\right\|_{L^2(X)\to L^2(X)}$$ $$\times\left\|\psi_{\delta'}(r)^{1/2}{\cal D}_r\left(\widetilde G-\lambda^2\mp 
i\varepsilon\right)^{-1}\psi_{\delta'}(r)^{1/2}
\right\|_{L^2(X)\to L^2(X)}$$ $$+O(\lambda)\sum_{\pm}\left\|\psi_{\delta'}(r)^{1/2}\left(\widetilde G-\lambda^2\pm 
i\varepsilon\right)^{-1}\psi_{\delta'}(r)^{1/2}
\right\|_{L^2(X)\to L^2(X)}$$ $$\times\left\|\psi_{\delta'}(r)^{1/2}r^{-1}\Lambda_w^{1/2}
\left(\widetilde G-\lambda^2\mp i\varepsilon\right)^{-1}\psi_{\delta'}(r)^{1/2}
\right\|_{L^2(X)\to L^2(X)},\eqno{(2.38)}$$
where we have used that $\psi_\delta\le\psi_{\delta'}$ for $\delta'\le\delta$. It is clear now that (1.11) 
with $\alpha_1=\alpha_2=0$ follows from (1.7) and (2.38). 
 Furthemore, it is easy to see that when $|\alpha_1|+|\alpha_2|\ge 1$ the estimate (1.11) follows from (1.7), 
(1.11) with $\alpha_1=\alpha_2=0$ and Lemma 2.1. Indeed, we have
$$\left\|\langle x\rangle^{-\frac{3+\delta'}{2}}\partial_x^{\alpha_1}\left(G-\lambda^2\pm i
\varepsilon\right)^{-2}\partial_x^{\alpha_2}
\langle x\rangle^{-\frac{3+\delta'}{2}}\right\|_{L^2\to L^2}$$ $$\le 
\left\|\langle x\rangle^{-\frac{3+\delta'}{2}}\partial_x^{\alpha_1}\left(G-i\lambda^2\right)^{-1}
\langle x\rangle^{\frac{3+\delta'}{2}}\right\|_{L^2\to L^2}$$ $$\times
\left\|\langle x\rangle^{-\frac{3+\delta'}{2}}\left(G-i\lambda^2\right)\left(G-\lambda^2\pm i\varepsilon\right)^{-2}
(G-i\lambda^2)\langle x\rangle^{-\frac{3+\delta'}{2}}\right\|_{L^2\to L^2}$$
 $$\times \left\|\langle x\rangle^{-\frac{3+\delta'}{2}}\partial_x^{\alpha_2}\left(G+i\lambda^2\right)^{-1}
\langle x\rangle^{\frac{3+\delta'}{2}}\right\|_{L^2\to L^2}$$
 $$\le O(\lambda^{|\alpha_1|+|\alpha_2|-4})\sum_{k=0}^2\lambda^{2k}\left\|\langle x\rangle^{-\frac{3+\delta'}{2}}
\left(G-\lambda^2\pm i\varepsilon\right)^{-k}
\langle x\rangle^{-\frac{3+\delta'}{2}}\right\|_{L^2\to L^2}$$ $$\le C\lambda^{|\alpha_1|+|\alpha_2|-2}.$$
\eproof

{\it Proof of Corollary 1.2.} 
We will use Theorem 1.1 and the observation that the constant $C$ in the right-hand
side of (1.7) depends only on the parameter $\delta'$, provided $\delta'\le \delta$ 
(an argument already used above in the case when $b^S\equiv V^S\equiv 0$ and which is true in the general case). 
Since by assumption $b(0)=0$ and the function $b$ is continuous in $r$, given any $\epsilon>0$ there is $0<\theta\le 1$
 so that $|b(x)|\le\epsilon$ for $|x|\le \theta$. Let $\zeta\in C_0^\infty({\bf R})$, $0\le\zeta\le 1$, 
$\zeta(\tau)=1$ for $|\tau|\le 1/2$, $\zeta(\tau)=0$ for $|\tau|\ge 1$. We are going to apply Theorem 1.1 to the operator
$$G_2=-\Delta+i(1-\zeta)(|x|/\theta)b(x)\cdot\nabla +i\nabla\cdot b(x)(1-\zeta)(|x|/\theta)+V(x)+\left|b(x)\right|^2.$$
Let $\chi\in C^\infty({\bf R})$, $0\le\chi\le 1$, $\chi(r)=0$ for $r\le r_0+1$, $\chi(r)=1$ for $r\ge r_0+2$. Set
$$\widetilde b^L(x)=\chi(|x|)b^L(x),\quad \widetilde V^L(x)=\chi(|x|)V^L(x),$$
$$\widetilde b^S(x)=(1-\zeta)(|x|/\theta)\left(b^S(x)+(1-\chi)(|x|)b^L(x)\right),$$
$$\widetilde V^S(x)=V^S(x)+(1-\chi)(|x|)V^L(x)+\zeta(|x|/\theta)(2-\zeta(|x|/\theta))|b(x)|^2.$$
It is easy to see that the operator $G_2$ is of the form
$$G_2=\left(i\nabla+\widetilde b^L+\widetilde b^S\right)^2+\widetilde V^L+\widetilde V^S,$$
and that the conditions of Corollary 1.2 imply that the functions $\widetilde b^L$, $\widetilde b^S$, $\widetilde V^L$
and $\widetilde V^S$ satisfy (1.1)-(1.6) with possibly a new constant $\delta>0$ independent of $\epsilon$. Therefore, by Theorem 1.1 the 
operator $G_2$ satisfies the estimate (1.7) with a constant $C$ in the
right-hand side independent of $\epsilon$. On the other hand, the difference $G-G_2$ is a first order differential operator
of the form $O(\epsilon)\cdot\nabla+\nabla\cdot O(\epsilon)$ with coefficients supported in $|x|\le 1$. 
Taking $\epsilon>0$ small enough, independent of $\lambda$, and proceeding in the same way as in the proof of (2.33)
above, we obtain that the operator $G$ satisfies (1.7), too. 
\eproof

{\it Proof of Corollary 1.3.} 
It is similar to the proof of Corollary 1.2 above. 
Since by assumption $b(x)=O(\langle x\rangle^{-\delta})$, given any $\epsilon>0$ there is $x_\epsilon\in{\bf R}^n$, 
$|x_\epsilon|\gg 1$, so that $|b(x)|\le\epsilon$ for $|x-x_\epsilon|\le 1$. We would like to apply Theorem 1.1 to the operator
$$G_3=-\Delta+i(1-\zeta)(|x-x_\epsilon|)b(x)\cdot\nabla +i\nabla\cdot b(x)(1-\zeta)(|x-x_\epsilon|)+V(x)+\left|b(x)\right|^2.$$
To this end, introduce the polar coordinates $r=|x-x_\epsilon|$, $w=\frac{x-x_\epsilon}{|x-x_\epsilon|}$. It is easy to see
that the conditions (1.14)-(1.18) imply that the coefficients of the operator $G_3$ satisfy (1.1)-(1.6) in these
new polar coordinates with the same constant $\delta>0$. Therefore, by Theorem 1.1 the operator $G_3$ satisfies the
estimate (1.7) with weights $\langle x-x_\epsilon\rangle^{-\frac{1+\delta'}{2}}$ and with a constant $C$ in the
right-hand side independent of $\epsilon$. On the other hand, the difference $G-G_3$ is a first order differential operator
of the form $O(\epsilon)\cdot\nabla+\nabla\cdot O(\epsilon)$ with coefficients supported in $|x-x_\epsilon|\le 1$. 
Taking $\epsilon>0$ small enough, independent of $\lambda$, and proceeding in the same way as in the proof of (2.33)
above, we obtain that the operator $G$ satisfies (1.7), too, with weights $\langle x-x_\epsilon\rangle^{-\frac{1+\delta'}{2}}
\sim \langle x\rangle^{-\frac{1+\delta'}{2}}$. 
\eproof

\section{Dispersive estimates}

Let $\varphi\in C_0^\infty((0,+\infty))$. It is easy to see that the estimates (1.23) and (1.24) follow from the following
semi-classical dispersive estimates (e.g. see Section 2 of \cite{kn:CV2}).

\begin{Theorem} Under the assumptions of Theorem 1.4, there exist  constants $C,h_0>0$ such that for all $0<h\le h_0$,
$t\neq 0$, we have the estimate
$$\left\|e^{it\sqrt{G}}\varphi(h\sqrt{G})\right\|_{L^1\to L^\infty}\le Ch^{-3}|t|^{-1}.\eqno{(3.1)}$$
Moreover, for every $\delta'>0$ there exist $C,h_0>0$ such that for all $0<h\le h_0$,
$t\neq 0$, we have the estimate
$$\left\|e^{it\sqrt{G}}\varphi(h\sqrt{G})\langle x\rangle^{-3/2-\delta'}\right\|_{L^2\to L^\infty}
\le Ch^{-2}|t|^{-1}.\eqno{(3.2)}$$
\end{Theorem}

{\it Proof.} We are going to use the formula
$$e^{it\sqrt{G_0}}\varphi(h\sqrt{G_0})=(\pi i)^{-1}\int_0^\infty e^{it\lambda}\varphi(h\lambda)
\left(R_0^+(\lambda)-R_0^-(\lambda)\right)\lambda
d\lambda,\eqno{(3.3)}$$
where $R_0^\pm(\lambda)=(G_0-\lambda^2\pm i0)^{-1}$ are the three dimensional outgoing and incoming free resolvents with kernels given by
$$[R_0^\pm(\lambda)](x,y)=\frac{e^{\pm i\lambda|x-y|}}{4\pi|x-y|}.$$
 We also have the formula
$$e^{it\sqrt{G}}\varphi(h\sqrt{G})=(\pi i)^{-1}\int_0^\infty e^{it\lambda}\varphi(h\lambda)
\left(R^+(\lambda)-R^-(\lambda)\right)\lambda
d\lambda,\eqno{(3.4)}$$
where $R^\pm(\lambda)=(G-\lambda^2\pm i0)^{-1}$ are the outgoing and incoming perturbed resolvents satisfying the relation
$$R^\pm(\lambda)-R_0^\pm(\lambda)=R_0^\pm(\lambda)LR^\pm(\lambda)=:T^\pm(\lambda)=T_1^\pm(\lambda)+T_2^\pm(\lambda),
\eqno{(3.5)}$$
where 
$$T_1^\pm(\lambda)=R_0^\pm(\lambda)L
R_0^\pm(\lambda),\quad T_2^\pm(\lambda)=R_0^\pm(\lambda)LR^\pm(\lambda)LR_0^\pm(\lambda),$$
$$L=G-G_0=ib(x)\cdot\nabla+i\nabla\cdot b(x)+|b(x)|^2+V(x).$$
In view of (3.3), (3.4) and (3.5) we can write
$$e^{it\sqrt{G}}\varphi(h\sqrt{G})-e^{it\sqrt{G_0}}\varphi(h\sqrt{G_0})=
(i\pi h)^{-1}\int_0^\infty e^{it\lambda}\widetilde\varphi(h\lambda)T(\lambda)d\lambda,\eqno{(3.6)}$$
where we have put $\widetilde\varphi(\lambda)=\lambda\varphi(\lambda)$, $T=T^+-T^-$. It is easy to see that the estimates (3.1) and (3.2)
follow from (3.6) and the following

\begin{prop} The operator-valued functions $T(\lambda):L^1\to L^\infty$ and $T(\lambda)\langle 
x\rangle^{-3/2-\delta'}:L^2\to L^\infty$ are $C^1$ for $\lambda$ large enough and satisfy the estimates (with $k=0,1$)
$$\left\|\partial_\lambda^k T(\lambda)\right\|_{L^1\to L^\infty}\le C\lambda,\eqno{(3.7)}$$
$$\left\|\partial_\lambda^k T(\lambda)\langle x\rangle^{-3/2-\delta'}\right\|_{L^2\to L^\infty}\le C.\eqno{(3.8)}$$
\end{prop}

{\it Proof.} We will need the following properties of the three dimensional free resolvent.

\begin{lemma} We have the estimates
$$\left\|\partial_\lambda^k\left(R_0^+(\lambda)-R_0^-(\lambda)\right)\right\|_{L^1\to L^\infty}\le C\lambda,\quad k=0,1,
\eqno{(3.9)}$$
$$\left\|\partial_\lambda^kR_0^\pm(\lambda)\langle x\rangle^{-1/2-k-\delta'}\right\|_{L^2\to L^\infty}
+\left\|\langle x\rangle^{-1/2-k-\delta'}\partial_\lambda^kR_0^\pm(\lambda)\right\|_{L^1\to L^2}\le C,\quad
k=0,1,\eqno{(3.10)}$$
$$\left\|\partial_x^\alpha\partial_\lambda R_0^\pm(\lambda)\langle x\rangle^{-3/2-\delta'}\right\|_{L^2\to L^\infty}
+\left\|\langle x\rangle^{-3/2-\delta'}\partial_\lambda R_0^\pm(\lambda)\partial_x^\alpha\right\|_{L^1\to L^2}\le C\lambda,
\quad |\alpha|=1.\eqno{(3.11)}$$
Moreover, if $|\alpha|=1$, given any $\gamma>0$ independent of $\lambda$ the operator $\partial_x^\alpha R_0^\pm(\lambda)$ 
can be decomposed as 
${\cal K}_{1,\alpha}^\pm(\lambda)+{\cal K}_{2,\alpha}$, where
$$\left\|{\cal K}_{1,\alpha}^\pm(\lambda)^*\langle x\rangle^{-1/2-\delta'}\right\|_{L^2\to L^\infty}
+\left\|\langle x\rangle^{-1/2-\delta'}{\cal K}_{1,\alpha}^\pm(\lambda)\right\|_{L^1\to L^2}\le C_\gamma\lambda,\eqno{(3.12)}$$
$$\left\|{\cal K}_{2,\alpha}^*\right\|_{L^\infty\to L^\infty}
+\left\|{\cal K}_{2,\alpha}\right\|_{L^1\to L^1}\le\gamma.\eqno{(3.13)}$$
\end{lemma}

{\it Proof.} The estimate (3.9) follows from the fact that the kernel of the operator $$\partial_\lambda^k\left(R_0^+(\lambda)-
R_0^-(\lambda)\right)$$ is $O(\lambda)$, while (3.10) follows from the fact that the kernel of the operator 
$\partial_\lambda^kR_0^\pm(\lambda)$
is $O\left(|x-y|^{k-1}\right)$ uniformly in $\lambda$. It is also easy to see that if $|\alpha|=1$, the kernel of
$\partial_x^\alpha\partial_\lambda R_0^\pm(\lambda)$ is $O(\lambda)$, which clearly implies (3.11). Furthemore, observe that
the kernel of $\partial_x^\alpha R_0^\pm(\lambda)$ is equal to
$$\frac{\partial^\alpha|x-y|}{\partial x^\alpha}\left(\pm i\lambda\frac{e^{\pm i\lambda|x-y|}}{|x-y|}-
\frac{e^{\pm i\lambda|x-y|}}{|x-y|^2}\right)$$
 $$=\frac{\partial^\alpha|x-y|}{\partial x^\alpha}\left(\pm i\lambda\frac{e^{\pm i\lambda|x-y|}}{|x-y|}+
\frac{1-e^{\pm i\lambda|x-y|}}{|x-y|^2}+\frac{\rho(|x-y|/\gamma')-1}{|x-y|^2}\right)$$ $$-
\frac{\partial^\alpha|x-y|}{\partial x^\alpha}
\frac{\rho(|x-y|/\gamma')}{|x-y|^2}:=K_{1,\alpha}^\pm(x,y)+K_{2,\alpha}(x,y),$$
where $\gamma'>0$ and $\rho\in C_0^\infty({\bf R})$, $0\le\rho\le 1$, 
$\rho(\sigma)=1$ for $|\sigma|\le 1$, $\rho(\sigma)=0$ for $|\sigma|\ge 2$. 
Denote by ${\cal K}_{1,\alpha}^\pm(\lambda)$ (resp. ${\cal K}_{2,\alpha}$) the operator with kernel $K_{1,\alpha}^\pm$
(resp. $K_{2,\alpha}$). Clearly, $K_{1,\alpha}^\pm=O_{\gamma'}(\lambda)|x-y|^{-1}$,
which implies (3.12). On the other hand, the left-hand side of (3.13) is upper bounded by
$$C\int_{{\bf R}^3} \frac{\rho(|x-y|/\gamma')}{|x-y|^2}dy\le C
\int_{|z|\le\gamma'}|z|^{-2}dz\le \widetilde C\gamma'.\eqno{(3.14)}$$
Choosing $\gamma'=\gamma/\widetilde C$ we get (3.13).
\eproof

Using Theorem 1.1 and Lemma 3.3 together with (1.20) and the fact that the operator $\partial_x^\alpha$ commutes
with the free resolvent, we obtain
$$\sum_{k=0}^1\left\|\partial_\lambda^kT(\lambda)\langle x\rangle^{-3/2-\delta'}\right\|_{L^2\to L^\infty}
\le \sum_\pm\left\|\frac{dR_0^\pm(\lambda)}{d\lambda}
LR^\pm(\lambda)\langle x\rangle^{-3/2-\delta'}\right\|_{L^2\to L^\infty}$$
 $$+\sum_{k=0}^1\left\|\left(R_0^+(\lambda)
L\frac{d^kR^+(\lambda)}{d\lambda^k}-R_0^-(\lambda)
L\frac{d^kR^-(\lambda)}{d\lambda^k}\right)\langle x\rangle^{-3/2-\delta'}\right\|_{L^2\to L^\infty}$$
$$\le C\sum_\pm\sum_{0\le|\alpha_1|+|\alpha_2|\le 1}\left\|\partial_x^{\alpha_1}\partial_\lambda R_0^\pm(\lambda)
\langle x\rangle^{-3/2-\delta/2}\right\|_{L^2\to L^\infty}\left\|\langle x\rangle^{-1/2-\delta/2}
\partial_x^{\alpha_2} R^\pm(\lambda)
\langle x\rangle^{-3/2-\delta'}\right\|_{L^2\to L^2}$$
 $$+ C\sum_\pm\sum_{k=0}^1\sum_{0\le|\alpha|\le 1}\left\|R_0^\pm(\lambda)
\langle x\rangle^{-1/2-\delta/2}\right\|_{L^2\to L^\infty}\left\|\langle x\rangle^{-3/2-\delta/2}
\partial_x^{\alpha}\partial_\lambda^k R^\pm(\lambda)
\langle x\rangle^{-3/2-\delta'}\right\|_{L^2\to L^2}$$
  $$+ C\sum_\pm\sum_{k=0}^1\sum_{|\alpha|= 1}\left\|{\cal K}_{1,\alpha}^\mp(\lambda)^*
\langle x\rangle^{-1/2-\delta/2}\right\|_{L^2\to L^\infty}\left\|\langle x\rangle^{-3/2-\delta/2}
\partial_\lambda^k R^\pm(\lambda)\langle x\rangle^{-3/2-\delta'}\right\|_{L^2\to L^2}$$ 
$$+ C\sum_{k=0}^1\sum_{|\alpha|= 1}\left\|{\cal K}_{2,\alpha}^*\right\|_{L^\infty\to L^\infty}
\left\|\left(\partial_\lambda^kR^+(\lambda)-
\partial_\lambda^k R^-(\lambda)\right)\langle x\rangle^{-3/2-\delta'}\right\|_{L^2\to L^\infty}$$
 $$\le C_\gamma+O(\gamma)\sum_{k=0}^1\left\|\partial_\lambda^kT(\lambda)
\langle x\rangle^{-3/2-\delta'}\right\|_{L^2\to L^\infty}\eqno{(3.15)}$$
for every $\gamma>0$. Taking $\gamma$ small enough we can absorb the second term in the right-hand side of (3.15)
and get (3.8). Let us see now that the operator $T_1=T_1^+-T_1^-$ satisfies (3.7). By Lemma 3.3 we have
$$\sum_{k=0}^1\left\|\partial_\lambda^kT_1(\lambda)\right\|_{L^1\to L^\infty}\le \sum_{k=0}^1\left\|\frac{d^kR_0^+(
\lambda)}{d\lambda^k}LR_0^+(\lambda)-\frac{d^kR_0^-(\lambda)}{d\lambda^k}LR_0^-(\lambda)
\right\|_{L^1\to L^\infty}$$ 
$$+\left\|R_0^+(\lambda)L\frac{dR_0^+(\lambda)}{d\lambda}-R_0^-(\lambda)L\frac{dR_0^-(\lambda)}
{d\lambda}\right\|_{L^1\to L^\infty}$$
 $$\le C\sum_\pm\sum_{k=0}^1\sum_{|\alpha|\le 1}\left\|\partial_x^\alpha\partial_\lambda^k R_0^\pm(\lambda)
\langle x\rangle^{-3/2-\delta/2}\right\|_{L^2\to L^\infty}\left\|\langle x\rangle^{-1/2-\delta/2}R_0^\pm(\lambda)
\right\|_{L^1\to L^2}$$ 
$$+C\sum_\pm\sum_{|\alpha|\le 1}\left\|R_0^\pm(\lambda)
\langle x\rangle^{-1/2-\delta/2}\right\|_{L^2\to L^\infty}\left\|\langle x\rangle^{-3/2-\delta/2}\partial_x^\alpha
\partial_\lambda R_0^\pm(\lambda)\right\|_{L^1\to L^2}$$
 $$+C\sum_\pm\sum_{k=0}^1\sum_{|\alpha|=1}\left\|\partial_\lambda^k R_0^\pm(\lambda)
\langle x\rangle^{-3/2-\delta/2}\right\|_{L^2\to L^\infty}\left\|\langle x\rangle^{-1/2-\delta/2}
{\cal K}_{1,\alpha}^\pm(\lambda)\right\|_{L^1\to L^2}$$ 
$$+C\sum_\pm\sum_{|\alpha|=1}\left\|{\cal K}_{1,\alpha}^\mp(\lambda)^*
\langle x\rangle^{-1/2-\delta/2}\right\|_{L^2\to L^\infty}\left\|\langle x\rangle^{-3/2-\delta/2}
\partial_\lambda R_0^\pm(\lambda)\right\|_{L^1\to L^2}$$
 $$+C\sum_{k=0}^1\sum_{|\alpha|=1}\left\|\partial_\lambda^k\left(R_0^+(\lambda)-R_0^-(\lambda)\right)\right\|_{L^1\to L^\infty}
\left(\left\|{\cal K}_{2,\alpha}\right\|_{L^1\to L^1}+\left\|{\cal K}_{2,\alpha}^*\right\|_{L^\infty\to L^\infty}\right)
\le C\lambda.\eqno{(3.16)}$$
Given a multi-index $\alpha=(\alpha_1,\alpha_2,\alpha_3)$ such that $|\alpha|\le 1$, define the function $b_\alpha$ as follows: $b_0=\left(|b|^2+V\right)/2$,
and if $|\alpha|=1$, $\alpha_j=1$, then $b_\alpha:=b_j$. 
The operator $T_2=T_2^+-T_2^-$ satisfies
$$\sum_{k=0}^1\left\|\partial_\lambda^kT_2(\lambda)\right\|_{L^1\to L^\infty}\le \sum_{k_1+k_2+k_3\le 1}
\left\|\sum_\pm \pm\frac{d^{k_1}R_0^\pm(\lambda)}{d\lambda^{k_1}}L\frac{d^{k_2}R^\pm(\lambda)}{d\lambda^{k_2}}L
\frac{d^{k_3}R_0^\pm(\lambda)}{d\lambda^{k_3}}\right\|_{L^1\to L^\infty}$$
 $$\le \sum_{k_1+k_2+k_3\le 1}\sum_{|\alpha_1|,|\alpha_2|,|\beta_1|,|\beta_2|\le 1,\,|\alpha_1|+|\alpha_2|\le 1,\,
|\beta_1|+|\beta_2|\le 1}A_{k_1,k_2,k_3}^{\alpha_1,\alpha_2,\beta_1,\beta_2}(\lambda)=:{\cal A}(\lambda),\eqno{(3.17)}$$
where 
$$A_{k_1,k_2,k_3}^{\alpha_1,\alpha_2,\beta_1,\beta_2}(\lambda)=\left\|\sum_\pm \pm\frac{d^{k_1}R_0^\pm(\lambda)}
{d\lambda^{k_1}}\partial_x^{\alpha_1}b_{\alpha_1,\alpha_2}\partial_x^{\alpha_2}\frac{d^{k_2}R^\pm(\lambda)}{d\lambda^{k_2}}
\partial_x^{\beta_2}b_{\beta_1,\beta_2}\partial_x^{\beta_1}
\frac{d^{k_3}R_0^\pm(\lambda)}{d\lambda^{k_3}}\right\|_{L^1\to L^\infty},$$
where $b_{\alpha_1,\alpha_2}=b_{\alpha_1}$ if $\alpha_2=0$, $b_{\alpha_1,\alpha_2}=b_{\alpha_2}$ if $\alpha_1=0$.
To bound these norms we will consider several cases.

Case 1. $\alpha_1=\beta_1=0$. By Theorem 1.1, Lemma 3.3 and (1.20), we have
$$A_{k_1,k_2,k_3}^{\alpha_1,\alpha_2,\beta_1,\beta_2}(\lambda)\le C\sum_\pm \left\|\frac{d^{k_1}R_0^\pm(\lambda)}
{d\lambda^{k_1}}\langle x\rangle^{-1/2-k_1-\delta/2}\right\|_{L^2\to L^\infty}$$
$$\times\left\|\langle x\rangle^{-1/2-k_2-\delta/2}\partial_x^{\alpha_2}\frac{d^{k_2}R^\pm(\lambda)}{d\lambda^{k_2}}
\partial_x^{\beta_2}\langle x\rangle^{-1/2-k_2-\delta/2}\right\|_{L^2\to L^2}\left\|\langle x\rangle^{-1/2-k_3-\delta/2}
\frac{d^{k_3}R_0^\pm(\lambda)}{d\lambda^{k_3}}\right\|_{L^1\to L^2}$$
 $$\le O(\lambda)\left\|\langle x\rangle^{-1/2-k_2-\delta/2}\partial_x^{\alpha_2}R^\pm(\lambda)^{1+k_2}
\partial_x^{\beta_2}\langle x\rangle^{-1/2-k_2-\delta/2}\right\|_{L^2\to L^2}\le C\lambda.\eqno{(3.18)}$$
 
Case 2. $|\alpha_1|+|\beta_1|\ge 1$, $k_1=1$ if $|\alpha_1|=1$ and $k_3=1$ if $|\beta_1|=1$. This case is treated in precisely
the same way as Case 1.

Case 3. $k_1=k_2=0$, $k_3=1$, $|\alpha_1|=1$, $\alpha_2=0$. By Theorem 1.1, Lemma 3.3 and (1.20), we have
$$A_{0,0,1}^{\alpha_1,0,\beta_1,\beta_2}(\lambda)\le C\sum_\pm \left\|{\cal K}_{1,\alpha_1}^\mp(\lambda)^*
\langle x\rangle^{-1/2-\delta/2}\right\|_{L^2\to L^\infty}$$
$$\times\left\|\langle x\rangle^{-1/2-\delta/2}R^\pm(\lambda)
\partial_x^{\beta_2}\langle x\rangle^{-1/2-\delta/2}\right\|_{L^2\to L^2}\left\|\langle x\rangle^{-3/2-\delta/2}
\partial_x^{\beta_1}\partial_\lambda R_0^\pm(\lambda)\right\|_{L^1\to L^2}$$
 $$+C\left\|{\cal K}_{2,\alpha_1}^*\right\|_{L^\infty\to L^\infty}\left\|\sum_\pm \pm R^\pm(\lambda)\partial_x^{\beta_2}
b_{\beta_1,\beta_2}\partial_x^{\beta_1}\partial_\lambda R_0^\pm(\lambda)\right\|_{L^1\to L^\infty}$$
 $$\le C_\gamma\lambda+O(\gamma)\left\|\sum_\pm \pm R_0^\pm(\lambda)\partial_x^{\beta_2}
b_{\beta_1,\beta_2}\partial_x^{\beta_1}\partial_\lambda R_0^\pm(\lambda)\right\|_{L^1\to L^\infty}$$
 $$+O(\gamma)\left\|\sum_\pm \pm R_0^\pm(\lambda)L R^\pm(\lambda)\partial_x^{\beta_2}
b_{\beta_1,\beta_2}\partial_x^{\beta_1}\partial_\lambda R_0^\pm(\lambda)\right\|_{L^1\to L^\infty}.\eqno{(3.19)}$$
In the same way as in the proof of (3.16) one can see that the second term in the right-hand side of (3.19) is
$O(\lambda)$. On the other hand, it is clear that the third one is bounded by $O(\gamma){\cal A}(\lambda)$. In other words, (3.19) yields
$$A_{0,0,1}^{\alpha_1,0,\beta_1,\beta_2}(\lambda)\le C_\gamma\lambda+O(\gamma){\cal A}(\lambda).\eqno{(3.20)}$$

Case 4. $k_1=1$, $k_2=k_3=0$, $|\beta_1|=1$, $\beta_2=0$. This case is treated in the same way as Case 3.

Case 5. $k_1=k_3=0$, $k_2=1$, $|\alpha_1|=|\beta_1|=1$, $\alpha_2=\beta_2=0$. By Theorem 1.1, Lemma 3.3 and (1.20), we have
$$A_{0,1,0}^{\alpha_1,0,\beta_1,0}(\lambda)\le C\sum_\pm \left\|{\cal K}_{1,\alpha_1}^\mp(\lambda)^*
\langle x\rangle^{-1/2-\delta/2}\right\|_{L^2\to L^\infty}$$
$$\times\left\|\langle x\rangle^{-3/2-\delta/2}\frac{dR^\pm(\lambda)}{d\lambda}
\langle x\rangle^{-3/2-\delta/2}\right\|_{L^2\to L^2}\left\|\langle x\rangle^{-1/2-\delta/2}
{\cal K}_{1,\beta_1}^\pm(\lambda)\right\|_{L^1\to L^2}$$
 $$+C\left\|{\cal K}_{2,\alpha_1}^*\right\|_{L^\infty\to L^\infty}\left\|\sum_\pm \pm \frac{dR^\pm(\lambda)}{d\lambda}
b_{\beta_1,0}\partial_x^{\beta_1}R_0^\pm(\lambda)\right\|_{L^1\to L^\infty}$$
 $$+C\left\|{\cal K}_{2,\beta_1}\right\|_{L^1\to L^1}\left\|\sum_\pm \pm R_0^\pm(\lambda)\partial_x^{\alpha_1}
b_{\alpha_1,0}\frac{dR^\pm(\lambda)}{d\lambda}\right\|_{L^1\to L^\infty}$$
 $$+C\left\|{\cal K}_{2,\alpha_1}^*\right\|_{L^\infty\to L^\infty}\left\|{\cal K}_{2,\beta_1}\right\|_{L^1\to L^1}
\left\|\sum_\pm \pm \frac{dR^\pm(\lambda)}{d\lambda}\right\|_{L^1\to L^\infty}.\eqno{(3.21)}$$
By (3.9), (3.16) and (3.17), we have
$$\left\|\frac{d(R^+(\lambda)-R^-(\lambda))}{d\lambda}\right\|_{L^1\to L^\infty}\le
 \left\|\frac{d(R_0^+(\lambda)-R_0^-(\lambda))}{d\lambda}\right\|_{L^1\to L^\infty}+
\left\|\frac{dT(\lambda)}{d\lambda}\right\|_{L^1\to L^\infty}$$
 $$\le C\lambda+\left\|\frac{dT_2(\lambda)}{d\lambda}\right\|_{L^1\to L^\infty}\le C\lambda+{\cal A}(\lambda).\eqno{(3.22)}$$
Similarly, one can easily see that the second and the third terms in the right-hand side of (3.21) are bounded
by $C\lambda+O(\gamma){\cal A}(\lambda)$. Thus we obtain
$$A_{0,1,0}^{\alpha_1,0,\beta_1,0}(\lambda)\le C_\gamma\lambda+O(\gamma){\cal A}(\lambda).\eqno{(3.23)}$$
Summing up the above inequalities we conclude
$${\cal A}(\lambda)\le C_\gamma\lambda+O(\gamma){\cal A}(\lambda).\eqno{(3.24)}$$
Taking $\gamma>0$ small enough, independent of $\lambda$, we can absorb the second term in the right-hand side of
(3.24) and conclude that ${\cal A}(\lambda)=O(\lambda)$. This together with (3.16) and (3.17) imply (3.7).
\eproof

{\bf Acknowledgements.} A part of this work has been carried out while F. C. and C. C. were visiting the Universit\'e de
Nantes, France. F. Cardoso has been partially supported by the agreement Brazil-France in Mathematics--Proc. 49.0733/2010-7,
and C. Cuevas has been partially supported by the CNRS-France. The first two authors are also partially supported by
the CNPq-Brazil.

F. Cardoso

Universidade Federal de Pernambuco, 

Departamento de Matem\'atica, 

CEP. 50540-740 Recife-Pe, Brazil,

e-mail: fernando@dmat.ufpe.br

\quad

C. Cuevas

Universidade Federal de Pernambuco, 

Departamento de Matem\'atica, 

CEP. 50540-740 Recife-Pe, Brazil,

e-mail: cch@dmat.ufpe.br

\quad

G. Vodev

Universit\'e de Nantes,

 D\'epartement de Math\'ematiques, UMR 6629 du CNRS,
 
 2, rue de la Houssini\`ere, BP 92208, 
 
 44332 Nantes Cedex 03, France,
 
 e-mail: vodev@math.univ-nantes.fr

\end{document}